\documentclass[a4paper,12pt,oneside]{article}
\usepackage[top=3cm, bottom=3cm, left=3cm, right=3cm,dvipdfm]{geometry}
\usepackage{t1enc}
\usepackage[latin2]{inputenc}
\usepackage{amsmath,amssymb,amsthm}
\usepackage[dvipdfm]{hyperref}
\usepackage{xcolor}

\hypersetup{
    colorlinks,
    linkcolor={blue!50!black},
    citecolor={red!40!black},
    urlcolor={green!60!black}
}
\usepackage{dsfont}
\usepackage{color}
\usepackage{graphicx}
\usepackage{enumerate}

\theoremstyle{plain}
\newtheorem{thm}{Theorem}[section]
\newtheorem{lemma}[thm]{Lemma}

\newtheorem{claim}[thm]{Claim}
\newtheorem{stm}[thm]{Statement}
\newtheorem{coroll}[thm]{Corollary}
\theoremstyle{definition}
\newtheorem{constr}{Construction}
\newtheorem{defi}[thm]{Definition}
\newtheorem{conj}{Conjecture}

\theoremstyle{remark}
\newtheorem{remark}[thm]{Remark}

\def\coex{\global\advance\count1by1}\count1=0
\def\cov{\mathrm{cov}}%
\def\COV{\mathrm{COV}}%
\def\Col{\mathrm{Col}}%
\def\KR{H_\mathrm{rich}}%
\def\prf{{\noindent \bf Proof.}\quad}%
\def\qed{\nobreak\hfill $\Box$\par\smallskip}%
\def\sqed{\nobreak\hfill {$\scriptstyle\Box$}}%
\def\ssqed{\nobreak\hfill {$\scriptscriptstyle\Box$}}%
\def\wl{wlog\ }%
\def\Wl{Wlog\ }%

\title{Covering complete partite hypergraphs by monochromatic components}

\author{Andr\'as Gy\'arf\'as\thanks{Research was supported in part by
grant  (no.\  K K104343) from the National Development Agency of
    Hungary, based on a source from the Research and Technology Innovation
    Fund.}\\[-0.8ex]
\small Alfr\'ed R\'enyi Institute of Mathematics\\[-0.8ex]
\small Hungarian Academy of Sciences\\[-0.8ex]
\small P.O. Box 127\\[-0.8ex]
\small Budapest, Hungary, H-1364\\[-0.8ex] \small
\texttt{gyarfas.andras@renyi.mta.hu} \and Zolt\'an Kir\'aly\thanks{Research
  was made during visiting Alfr\'ed R\'enyi
  Institute of Mathematics, Hungarian Academy of Sciences
  and was partially supported by grant
  (no.\  K 109240) from the National Development Agency of
    Hungary, based on a source from the Research and Technology Innovation
    Fund.} \\[-0.8ex]
\small E\"otv\"os Lor\'and University \\[-0.8ex]
\small Department of Computer Science  and\\[-0.8ex]
\small MTA-ELTE Egerv\'ary Research Group\\[-0.8ex]
\small P\'azm\'any P\'eter s\'et\'any 1/C \\[-0.8ex]
\small Budapest, Hungary, H-1117\\[-0.8ex]
\small \texttt{kiraly@cs.elte.hu}\\
}

\begin{document}
\thispagestyle{plain}
\maketitle

\begin{abstract} A well-known special case of a conjecture attributed to Ryser
  (actually 
appeared in the thesis of Henderson \cite{HE}) states that $k$-partite
intersecting hypergraphs have transversals of at most $k-1$ vertices. An
equivalent form of the conjecture in terms of coloring of complete graphs
is formulated in \cite{GY1}: if the edges of a complete graph $K$ are
colored with $k$ colors then the vertex set of $K$ can be covered by at
most $k-1$ sets, each connected in some color. It turned out that the
analogue of the conjecture for hypergraphs can be answered: Z. Kir\'aly
proved \cite{KI} that in every $k$-coloring of the edges of the $r$-uniform
complete hypergraph $K^r$ ($r\ge 3$), the vertex set of $K^r$  can be
covered by at most $\lceil k/r \rceil$ sets,   each connected in some
color.

Here we investigate the analogue problem for complete $r$-uniform
$r$-partite hypergraphs. An edge coloring of a hypergraph is called {\bf
spanning} if every vertex is incident to edges of any color used in the
coloring. We propose the following analogue of  Ryser conjecture.
\bigskip

\noindent {\em In every spanning $(r\!+\!t)$-coloring of the edges of a
  complete $r$-uniform
  $r$-partite hypergraph, the vertex set can be covered by at most
$t+1$ sets, each connected in some color.}

\bigskip

We show that the conjecture (if true) is best possible.
Our main result is that the conjecture is true
for $1\le t\le r-1$. We also prove a slightly weaker result for
$t\ge r$, namely that $t+2$ sets, each connected in some color, are enough to cover the vertex set.

To build a bridge between complete $r$-uniform and complete
$r$-uniform $r$-partite hypergraphs, we introduce a new notion.
A hypergraph is complete $r$-uniform
$(r,\ell)$-partite if it has all $r$-sets that intersect
each partite class in at most $\ell$ vertices (where $1\le\ell\le r$).

Extending our results achieved for $\ell=1$, we prove that
for any $r\ge 3,\; 2\le\ell\le r,\; k\ge 1+r-\ell$, in every spanning
$k$-coloring of the edges of a complete $r$-uniform
$(r,\ell)$-partite hypergraph, the vertex set can be
covered by at most $1+\lfloor \frac{k-r+\ell-1}{\ell}\rfloor$ sets, each
connected in some color.

\end{abstract}

\section{Introduction}

For an edge-colored hypergraph $H$ let $H_i$ denote its subhypergraph
consisting of edges colored by $i$. The connected components of $H_i$
are called monochromatic components of color $i$, and
a {\bf monochromatic component} refers to a monochromatic component of color
$i$ for some $i$. Here
 connectivity is understood in its weakest sense, a hypergraph is connected if
 either it has only one vertex or any two distinct vertices can be connected
 by a sequence of edges each intersecting the next. Every hypergraph can be
 uniquely partitioned into connected components. Components with a single
 vertex are called {\em trivial}.

 Given an edge-colored hypergraph $H$, let $c(H)$ denote the minimum integer
 $m$ such that $V=V(H)$, the vertex set of $H$, can be covered by $m$
 monochromatic components of $H$. An edge coloring of a hypergraph is called
 {\bf spanning} if every vertex is incident to edges of any color used in the
 coloring. Note that in spanning colorings every monochromatic component is
 non-trivial. The importance of this definition is shown in Theorem
 \ref{nonspanning}.

A conjecture attributed to Ryser which actually appeared in \cite{HE} is that
$k$-partite intersecting hypergraphs have  transversals of at most $k-1$
vertices.  An equivalent form is formulated in \cite{GY1} as follows: if $K$
is a complete graph with a $k$ coloring on its edges, then $c(K)\le k-1$. The
conjecture is true for $k\le 5$ and seems very difficult in general (further
information can be found in
\cite{EGYP}, \cite{GY2}). A particular feature of the conjecture is that
$c(K)\le k$ is obvious since the monochromatic stars at
any vertex form monochromatic components. Note that the conjecture is obvious
for colorings that are not spanning.

Surprisingly, the problem for hypergraphs is easier, Z. Kir\'aly in \cite{KI}
showed that
if the edges of the complete $r$-uniform hypergraph
$K$ ($r\ge 3$) are colored with $k$ colors, then $c(K)\le \lceil k/r \rceil$
and this is best possible (the $k=r$ case were already in \cite{GY1} extending
the well-known remark of Erd\H os and Rado stating that a graph or its
complement is connected).

The problem naturally extends for sparser host graphs (or
hypergraphs). Gy\'arf\'as and Lehel conjectured that for
$k$-colored complete bipartite graphs $G$, $c(G)\le 2k-2$ (see \cite{CFGYLT}),
here again $c(G)\le 2k-1$ is obvious.
For the hypergraph case \cite{FFGYT1,FFGYT2} initiated the study of $c(H)$
when $H$ has bounded independence number.

The main subject of the present paper is the case when the target hypergraph
$K$ is a complete $r$-uniform
$r$-partite hypergraph, i.e., when $V=V(K)$ is partitioned
into nonempty classes $V_1\cup\ldots\cup V_r$ and the edges of $K$ are the
sets containing one vertex from each class.  Let $\cov(r,k)$ denote the
maximum of $c(K)$ when $K$ ranges over spanning $k$-colorings of complete
$r$-uniform $r$-partite hypergraphs, and $\COV(r,k)$ denote the
maximum of $c(K)$ when $K$ ranges over (not necessarily spanning)
$k$-colorings of complete $r$-uniform $r$-partite hypergraphs.

Throughout the paper we {\bf always assume ${r\ge 3}$}. Our introductory
theorem shows that only the spanning colorings are the interesting ones.
For any positive integer $k$ we use the standard notation
$[k]=\{1,2,\ldots,k\}$.

\begin{thm}\label{nonspanning}
  If $r\ge 3$, then $\COV(r,k)= k$.
\end{thm}

\prf Let $K$ be a $k$-edge-colored $r$-uniform $r$-partite complete hypergraph.
Take an edge $e$ of $K$. Let $C_1,\ldots,C_\ell$ be the monochromatic
components with $|C_i\cap e|\ge r-1$. As $r>2$, clearly no two of them have the
same color, so $\ell\le k$. For every vertex $v\in V$ there is an edge $f\ni v$
with $|f\cap e|=r-1$, so $v$ is covered by one of these components.\\
For the sharpness let $V_1=[k]$ and color each edge $e$ by color $e\cap V_1$.
\qed

\medskip

We remark that if a coloring of the  $r$-uniform $r$-partite complete
hypergraph is spanning, then
\emph{all monochromatic components meet every class}.
An edge of color $i$ in a $k$-colored $r$-uniform hypergraph $K$
is called {\bf essential} if it is not contained in monochromatic components
of any color different from $i$.
When $\cov(r,k)$ is studied we may restrict ourselves to colorings having at
least one essential edge in every used color, since
otherwise a color can be eliminated by recoloring all edges of that color to
some other color and the resulting hypergraph
would  still have a spanning coloring and the same set of (maximal) monochromatic
components. This concept is established in \cite{KI}
and works well in the proof of our initial result.

\begin{thm}\label{lemma_essential}
$\cov(r,k)= 1$ for every $1\le k\le r\ge 3$.
\end{thm}

\prf
Let $e=\{v_1,\ldots,v_r\}$ be an essential edge of color $1$ in a
complete $r$-uniform $r$-partite
hypergraph  with vertex set $V=\cup _{i=1}^r V_i$ where $v_i\in V_i$.
Let $R_i=e-\{v_i\}$ and denote by $\Col(R_i)\subseteq [k]$
the set of colors appearing on any
edge of the form $R_i\cup\{v'_i\}$ (where $v'_i\in V_i$).
As $\Col(R_i)\cap \Col(R_j)=\{1\}$ for $i\ne j$, by the pigeonhole principle
there exists  $j$ such that  $\Col(R_j)=\{1\}$. Now $V_j$ is covered by the
monochromatic component containing $e$ (of color 1), and, as the coloring is
spanning, it necessarily covers the whole $V$. \qed

By Theorem \ref{lemma_essential} from this point
we may assume that $k=r+t$ with some integer $t\ge 1$.

\begin{conj}\label{conj}
$\cov(r,r+t)= t+1$ for every $r\ge 3,\; t\ge 1$.
\end{conj}

It is worth formulating this conjecture in dual form. Assume $K$ is  a
complete $r$-uniform $r$-partite hypergraph with a spanning
$k$-coloring. Consider a
new hypergraph $H$ with vertex set $V(K)$ whose edges are the vertex sets of the
monochromatic components in the coloring. The dual $F$ of this new
hypergraph $H$
is a $k$-uniform $k$-partite complete hypergraph whose edges are
partitioned
into $r$ classes with the property that any $r$ edges from different partite
classes have nonempty intersection. As the coloring  of $K$ was spanning,
monochromatic components have at least $r$ vertices.
In this setting Conjecture \ref{conj} can
be stated in terms of the transversal number $\tau(F)$, the minimum number of
vertices intersecting all edges of $F$.

\begin{conj}\label{conj1}
Assume that the edges of a $k$-uniform $k$-partite hypergraph $F$
with minimum degree at least $r\ge 3$  are
partitioned into $r$ classes so that any $r$ edges from different classes
have nonempty intersection.  Then $\tau(F)\le k-r+1$.
\end{conj}

In Section \ref{comppartite} we show that Conjecture \ref{conj} (if true) is
best possible, and it is ``almost'' true, i.e.,
$\cov(r,r+t)\le t+2$ for every $t\ge 1$ (Theorem \ref{large_t}).
We also prove that the conjecture is true for $1\le t \le r-2$ (Theorem
\ref{small_t}). Our most difficult result makes one further step, proving
Conjecture \ref{conj} for $t=r-1$ (Theorem \ref{3plus2}).

In Section \ref{trans} we investigate $c(H)$ for hypergraphs ``between''
complete and complete partite, in order to build a bridge between the results
proved in Section \ref{comppartite} and the results of \cite{KI}.
We call a hypergraph {\em $(r,\ell)$-partite} if its vertex set is partitioned
into $r$ nonempty classes, such that the intersection of any edge and any class
has at most $\ell$ vertices.
We call a hypergraph {\em complete $r$-uniform  $(r,\ell)$-partite} if it
contains all $r$-element sets as
edges which meet every partition class in at most $\ell$ vertices.
Let $\cov(r,\ell,k)$ denote the minimum number of
monochromatic components
needed to cover the vertex set of any complete $r$-uniform  $(r,\ell)$-partite
hypergraph in any spanning $k$-coloring.
For $2\le \ell\le r$ we determine exactly the values of $\cov(r,\ell,k)$.
We conclude our paper by summarizing the results achieved. Our main result is Theorem
\ref{r_l_e}, stating that
  $$\cov(r,\ell,k)= 1+\Bigl\lfloor \frac{k-r+\ell-1}{\ell}\Bigr\rfloor$$
for every $r\ge 3,\; k\ge 1+r-\ell, \; 1\le \ell\le r$,
\emph{except} for the cases
$(\ell=1$ and $k\ge 2r)$, where only we could prove a slightly
weaker upper bound.

\section{Results for complete  \texorpdfstring{$r$}{r}-uniform
  \texorpdfstring{$r$}{r}-partite
    hypergraphs}\label{comppartite}

\subsection{Lower bound}

\begin{constr}\label{basic_constr}
For $t\ge 1,\; r\ge 3,\; k=r+t$, we
define a complete $r$-uniform  $r$-partite hypergraph
$K(r,t)$ with a $k$-coloring of its edges as follows.  The vertex set $V$ of
$K(r,t)$ is partitioned into $r$ classes, $V_1,\dots,V_r$.
The first class $V_1$ has ${k\choose t}$ vertices associated to the
$t$-element subsets of $[k]$.  For $2\le j \le r$ set
$V_j=A_j^1\cup \dots \cup A_j^{k}$,
where the $A_j^i$-s are disjoint and have ${k-1\choose t-1}$
vertices. Fix an arbitrary linear order on every $A_j^i$.

First we define special edges of color $i$ for any $i\in [k]$. Consider the
set $W_i$ of ${k-1\choose  t-1}$
vertices of $V_1$ associated to $t$-sets of $[k]$ containing $i$.
\begin{itemize}
\item {\em Special edges of color $i$} are the
${k-1\choose  t-1}$ edges whose vertex from $W_i$ is the $\ell$-th in
  lexicographic order,
  and for all $2\le j\le r$ whose vertex from $V_j$ is the $\ell$-th
in the fixed linear order of $A_j^i$ for $\ell=1,\dots,{k-1\choose t-1}$.
Thus special edges of color $i$ form a matching for all $i$, $i=1,\dots,k$.
\item {\em Non-special edges} with vertices $v_1\in V_1,\dots,v_r\in V_r$ get
  their color as
the smallest $c\in [k]$ such that $c$ is not in the set associated to $v_1$
and $v_j\notin A_j^c$ for all $2\le j \le r$.
\end{itemize}
Note that every non-special
$r$-tuple $v_1,\dots,v_r$ gets a color because the conditions forbid at most
$t+r-1$ colors. Observe also that a special edge of color $i$ is always
disjoint from any other edge of color $i$. Consequently a special edge of
color $i$ forms a monochromatic component of color $i$ having $r$ vertices,
we call them small monochromatic components.

We claim that the coloring given is spanning. Suppose first that $v\in
V_1$ representing \wl the set $[t]\subset [k]$.  For any $1\le i \le t$,  $v$
is in a special edge of color $i$. On the other hand, for any $t<i\le r+t$ we
can select vertices  $v_2\in A_2^{j_2}\dots,v_r\in  A_r^{j_r}$ so that the
upper indices $j_t$ take all values except $i$ from $t+1,\dots,t+r$. Then the
non-special edge $v,v_2,\dots,v_r$ is colored by $i$.

On the other hand, let $v\in A_j^i$ for some $1<j\le r,\; 1\le i \le k$. Clearly
$v$ is in a special edge of color $i$.
For any $c\ne i$ such that $1\le c \le k$ we can take any vertex $w\in V_1$
associated to a $t$-set $A$ of
$[k]$ such that $c,i\notin A$. Set $B=[k]\setminus (A\cup \{i\} \cup \{c\})$.
Then from the $(r\!-\!2)$ $V_t$-s
where $t\notin \{1,j\}$ we can pick a set of $r-2$ vertices with distinct
superscripts in $B$.
These vertices together with $v,w$ define an edge that must be colored with
$c$. Thus the coloring of $K(r,t)$ is spanning.
\end{constr}

\begin{thm}\label{constr}
  $\cov(r,r+t)\ge t+1$ for every $r\ge 3,\; t\ge 1$.
\end{thm}

\prf Consider the hypergraph $K(r,t)$.
Note that the union of at most $t$ large monochromatic components
do not cover $V_1$. Let their colors be $c_1,\ldots,c_s$ with $s\le t$, and
take any $t$-set that contains $\{c_1,\ldots,c_s\}$; the vertex in $V_1$
associated to this set is not covered.

The uncovered vertices of $V_1$
must be covered by small monochromatic components, and every such component
can contain just one vertex of $V_1$.
Therefore
we need ${k-s\choose t-s}>t-s$ small monochromatic components to cover
them, thus altogether we need more than $s+(t-s)=t$ monochromatic components
to cover the vertices of $K(r,t)$. \qed

\medskip
\subsection{Upper bounds}

We need some additional notation.
We assign vectors of length $k$ to every element of the base
set $V=V_1\cup\ldots\cup V_r$. For $v\in V$ the $i$th coordinate
${\mathbf  v}(i)$ of
the associated vector
${\mathbf v}$ is the serial number of the monochromatic component of color $i$
containing $v$. The Hamming distance of two vertices
$\delta(v,w)=\delta({\mathbf v},{\mathbf w})$ is the number of places the two associated vectors differ.

\begin{stm}\label{stm}
  For $i=1,\ldots,r$ let $v_i\in V_i$. Then there exists  $c\in[k]$ and
  an integer $s$, such that ${\mathbf v}_i(c)=s$ for all $i\le r$.
\end{stm}

\prf
The edge $e=\{v_1,v_2,\ldots,v_r\}$ is colored by a color, say, by color $c$.
Then the vertices of $e$ belong to the same monochromatic component of color
$c$. \sqed

\begin{lemma}\label{base_lemma}
  Either $\cov(r,r+t)= 1$, or for any two vertices $v,w$ from different classes, 
  $\delta(v,w)\le t+1$.
\end{lemma}

\prf
\Wl $v\in V_1,\; w\in V_2$ and ${\mathbf  v}=1\ldots1$ and
${\mathbf  w}=1\ldots12\ldots2$, where the number of ones is at
most $r-2$. As the coloring is spanning and no monochromatic component covers
$V$, we can choose $v_3,\ldots,v_r$, such that $v_i\in V_i$ and
${\mathbf  v}_i(i-2)>1$.
However, this contradicts to Statement \ref{stm}. Thus the number of twos in
${\mathbf w}$ is at most $t+1$, so $\delta(v,w)\le t+1$. \qed

\begin{lemma}\label{gen_base_lemma}
If $\cov(r,r+t)> 1$ and
$w_1,\ldots,w_\ell$ are vertices from different classes, then for
$J=\{j\in[k]\;|\;{\mathbf w}_1(j)={\mathbf w}_2(j)=\ldots={\mathbf w}_\ell(j)\}$
we have $|J|\ge r+1-\ell$.
\end{lemma}

\prf
If $\ell=r$, then this statement coincides with Statement \ref{stm}.
Otherwise suppose $|J|\le r-\ell$ and $J=\{j_1,\ldots,j_{|J|}\}$.
We may choose at most $r-\ell$ vertices $u_1,\ldots,u_{|J|}$ from the
classes not having a $w_i$ with ${\mathbf u}_i(j_i)\ne {\mathbf w}_1(j_i)$,
contradicting to Statement \ref{stm}.
\qed

\begin{thm}\label{small_t}
  $\cov(r,r+t)\le t+1$ for every $1\le t\le r-2$ and $r\ge 3$.
\end{thm}

\prf
Suppose the statement does not hold.
First we claim that for any $i\ne j$ and for any $a\in V_i,\; b\in V_j$
we have $\delta(a,b)\le t$. Suppose not, \wl $a\in V_1,\; b\in V_2$, such that
${\mathbf  a}=1\ldots1$ and ${\mathbf  b}=1\ldots12\ldots 2$, where ${\mathbf b}$ ends
with $q$ 2-values, and  $q\ge t+1$, consequently $q=t+1$ by Lemma
\ref{base_lemma}.

As two monochromatic components do not cover $V$, there exists $d\in V_3$,
such that ${\mathbf d}(t+r)>2$.
By the assumption we have a vertex $c\in V_i$
for some $i$ with ${\mathbf c}(j)\ne 1$ for $j=1,\ldots,t+1$.

If $i>1$, then by Lemma \ref{base_lemma} $\delta(a,c)\le t+1$ and so
${\mathbf c}(j)= 1$ for $j=t+2,\ldots,t+r$. As $t+1\le r-1$,
$\delta(b,c)\ge 2t+2$, so $i=2$.
Now $\delta(d,b)\le t+1$ and $\delta(d,c)\le t+1$ but $\delta(b,c)=2t+2$, so
${\mathbf d}$ has to agree with either ${\mathbf b}$ or
${\mathbf c}$ in every coordinate where ${\mathbf b}$ and ${\mathbf c}$
differ. However, this is not the
case for the $(t+r)$-th coordinate.

If $i=1$, then by Lemma \ref{base_lemma} $\delta(b,c)\le t+1$ and so
${\mathbf c}(j)= {\mathbf b}(j)$ for $j=t+2,\ldots,t+r$, as $t+1\le r-1$.
Now $\delta(d,a)\le t+1$ and $\delta(d,c)\le t+1$ but $\delta(a,c)=2t+2$, so
${\mathbf d}$ has to agree with either ${\mathbf a}$ or
${\mathbf c}$ in every coordinate where ${\mathbf a}$ and ${\mathbf c}$
differ. However, this is not the
case for the $(t+r)$-th coordinate and the claim is proved.

Let $a,b,c$ as before, now $q\le t$, so we have
$\delta(a,c)\ge t+1$ and $\delta(b,c)\ge t+1$ but this contradicts to the claim
because either $c\not\in V_1$ or $c\not\in V_2$. \qed

\begin{thm}\label{large_t}
  $\cov(r,r+t)\le t+2$ for every $2\le r-1\le t$.
\end{thm}

\prf
Suppose the statement does not hold.
First we claim that there exist $i\ne j,\; a\in V_i,\; b\in V_j$, such that
$\delta(a,b)\ge r-2$. \Wl $a_\emptyset\in V_1$ with
${\mathbf  a}_\emptyset=11\ldots1$. For each $J\subseteq[k],\; |J|=t+2$,
there exists a vertex $a_J$ with ${\mathbf  a}_J(j)\ne 1$ for each $j\in J$.
(Note, that for $J\ne J'$, $a_J=a_{J'}$ is possible.)
These vertices have Hamming distance $\delta(a_J,a_\emptyset)>t+1$ from
$a_\emptyset$,
consequently, by Lemma \ref{base_lemma}, they all are in $V_1$.
Take any vertex $b\in V_2$, we claim that $\delta(b,a_J)\ge r-2$ for either
$J=\emptyset$ or for a $|J|=t+2$. Let $I=\{i\; | \; {\mathbf  b}(i)= 1\}$.
If $|I|<t+2$, then  $\delta(b,a_\emptyset)\ge r-2$ and we are done, otherwise take
any $J\subseteq I$ with $|J|=t+2$. Now obviously $\delta(b,a_J)\ge t+2\ge r+1$,
proving the claim.

By the claim  we have \wl $b_1\in V_1,\; b_2\in V_2$ where ${\mathbf
  b}_1=11\ldots1$ and
${\mathbf  b}_2=22\ldots211\ldots1$, and ${\mathbf  b}_2$ starts with $r-2+q$
twos ($q\ge 0$). If two monochromatic components cover $V$, then we are
done, otherwise we have $b_i\in V_i$ for $i=3,4,\ldots,r$, such that
${\mathbf  b}_i(i-2)>2$. Take also a vertex $d\in V$, where
${\mathbf  d}(j)\ne {\mathbf  b}_r(j)$ for $r-2< j\le r-2+q$ and
${\mathbf  d}(j)\ne 1$ for $r-2+q< j\le r+t$; this involves
$(r+t)-(r-1)+1=t+2$ coordinates, by our assumption such a vertex
must exist.

Thus $d\in V_i$ for some $1\le i\le r$.
Take the edge $e=\{b_j\; | \; j\ne i\}\cup\{d\}$, it is colored by some color
$c$.
Observe that ${\mathbf  d}$ differs form
both ${\mathbf  b}_1$ an ${\mathbf  b}_2$ in the last $t+2-q$ coordinates, so
$c\le r-2+q$. If $c\le r-2$, then we have $i_1,i_2\in [r]-\{i\}$ such that
${\mathbf  b}_{i_1}(c)\ne {\mathbf  b}_{i_2}(c)$, so
$r-2< c \le r-2+q$. However, if $i<r$, then
${\mathbf  b}_{r}(c)\ne {\mathbf  d}(c)$, otherwise
${\mathbf  b}_{1}(c)\ne {\mathbf  b}_{2}(c)$. Thus $c$ does not exist,
contradiction. \qed

\subsection{The case \texorpdfstring{$t=r-1$}{t=r-1}}

\begin{thm}\label{3plus2}  $\cov(r,2r-1)=r$ if $r\ge 3$.
\end{thm}

Suppose the statement does not hold, let $k=2r-1$ and
fix a $k$-colored $r$-uniform  $r$-partite hypergraph $K$
where $c(K)\ge r+1$ (and the coloring is spanning).

\begin{claim}\label{dist23}
  For any $i\ne j$ and $a\in V_i, b\in V_j$ we have
  $$ r-1\le \delta(a,b)\le r.$$
\end{claim}
\prf The upper bound comes from Lemma \ref{base_lemma}. To prove the lower
bound, \wl assume that we have $a\in V_1,\; b\in V_2$ with vectors
$${\mathbf a}=11\dots 111, \; {\mathbf b}=2\ldots21\ldots1$$
where ${\mathbf b}$ begins with $q\le r-2$ twos.
Suppose first that $q>0$.

As two monochromatic components do not cover $V$,
we can choose vertices $d_3,\ldots,d_r$ with $d_i\in V_i$ and
${\mathbf d}_i(i-2)=3$.
We claim that ${\mathbf d}_i(j)=1$ for all
$3\le i\le r$ and $r-1\le j\le 2r-1$. Otherwise the
index set
$J=\{j\;|\;  {\mathbf d}_3(j)={\mathbf d}_4(j)=\ldots={\mathbf d}_r(j)=1\}$
has size at most $r$, so there is a set $I$ such that
$J\subseteq I\subseteq \{r-1,\ldots,2r-1\}$ and $|I|=r$.
There is a vertex $c_I$ with the property
${\mathbf c}_I(j)\ne 1$ for all $j\in I$. As either $\delta(c_I,a)>r$ or
$\delta(c_I,b)>r$, $c_I\in V_1\cup V_2$. If $c_I\in V_1$, then
$c_I,b,d_3,\ldots,d_r$, otherwise $a,c_I,d_3,\ldots,d_r$ contradicts to
Lemma \ref{base_lemma}, so ${\mathbf d}_3(j)=1$ for all $r-1\le j\le 2r-1$.
Now for $I=\{r,\ldots,2r-1\}$ we also have  $c_I\in V_1\cup V_2$, and
${\mathbf c}_I(1)={\mathbf a}(1)=1$ (if $c_I\in V_2$) or
${\mathbf c}_I(1)={\mathbf b}(1)=2$ (if $c_I\in V_1$) also follows, so
$\delta(c_I,d_3)\ge 1+|I|=1+r$ contradicting to Lemma \ref{base_lemma}.

We conclude that $q=0$, thus
${\mathbf a}={\mathbf  b}$ are both the all-1 vectors.
Then for all $I\subset [2r-1],\; |I|=r$ there exist vertices $c_I$ such that
${\mathbf c}_I(j)\ne 1$ for all $j\in I$. As either $\delta(c_I,a)\le r$ or
$\delta(c_I,b)\le r$ must hold, ${\mathbf c}_I(j)= 1$ for all $j\in [2r-1]-I$.
Suppose that for $I_1,I_2\subset [2r-1]$ the complementary sets
$\overline{I_1},\overline{I_2}$ are disjoint. Then the corresponding vertices
$c_{I_1}$ and $c_{I_2}$
must be in the same vertex class, otherwise the Hamming distance of
their vectors would be at least $2(r-1)\ge r+1$. Since the Kneser graph defined
by disjoint $(r\!-\!1)$-element subsets of a $(2r\!-\!1)$-element ground set is
a connected graph, all the $c_I$-s are in the same class, call it the
full class; by symmetry we may assume that it is not $V_1$.
Select vertices $d_2\in V_2$ and $d_3\in V_3$  such that ${\mathbf d}_i(i)=3$.
Observe that ${\mathbf d}_i$ has at most $r$ ones
because $1\le\delta(a,d_i)\le r-2$
cannot happen.
If the full class is $V_2$, let $I$ contain the positions where
${\mathbf d}_3$ is 1, then
$a,c_I,d_3$ violate Statement \ref{stm}.
If the full class is $V_j$ for $j\ge 3$, let $I$ contain the positions where
${\mathbf d}_2$ is 1, now
$a,d_2,c_I$ violate Statement \ref{stm}.
%Ezt meg szebben kellene!
\sqed

\begin{claim}\label{dist2comp}
  If $a,b$ are two vertices from different partite
  classes such that $\delta(a,b)=r-1$, then some of these
  classes contain two vertices with Hamming distance $2r-1$.
\end{claim}

\prf Assume \wl that there are $a\in V_1,\; b\in V_2$ such that
${\mathbf a}=11\ldots 111,\; {\mathbf b}=2\ldots 21\ldots 1$
where ${\mathbf b}$ ends with exactly $r$ ones. There exists a
vertex $c\in V$ with ${\mathbf c}(j)\ne 1$ if $r\le j\le 2r-1$.
By Statement \ref{stm}, $c\in V_1\cup V_2$.
If $c\in V_1$, then (as $\delta(c,b)\le r$)
its vector starts with $r\!-\!1$ twos, so
$a,c$ is a pair required.
If $c\in V_2$, then its vector starts with $r\!-\!1$ ones, so
$b,c$ is a pair required.
\sqed

\begin{claim}\label{nocomp}
  For any two vertices $v,w$ from the same partite class,
  $$\delta({\mathbf v},{\mathbf w})<2r-1.$$
\end{claim}

\prf Assume indirectly that we have two vertices  \wl $v,w\in V_1$,
$${\mathbf  v}=11\ldots 11, \; {\mathbf w}=22\ldots 22.$$
For any $2\le i\le r, \; 1\le j\le 2r-1$ there exist vertices
$v_i^j\in V_i$ such that ${\mathbf   v}_i^j(j)=3$
by our assumption. They are all distinct because
their vectors must contain exactly $r\!-\!1$ ones and exactly
$r\!-\!1$ twos since their
distance from both $v,w$ must be at most $r\!-\!1$.

\begin{stm}\label{parity}
  Let ${\mathbf v}$ and ${\mathbf w}$ be $1$-$2$ vectors of the same length.
  If the number of ones in  ${\mathbf v}$ and ${\mathbf w}$ have the same
  parity, then $\delta({\mathbf v},{\mathbf w})$ is even, otherwise it is odd.
  \ssqed
\end{stm}

\begin{stm}\label{notcompl_in_same_matrix}
$\delta(v_i^j,v_{i}^{j'})\le 2r-2$ for any
$2\le i \le r, \; 1\le j<j'\le 2r-1$.
\end{stm}
\prf Suppose \wl  $\delta(v_2^1,v_{2}^{2})= 2r-1$ and
$${\mathbf v}_2^1=31\ldots 12\ldots 2,\;
     {\mathbf v}_{2}^{2}=232\ldots 21\ldots 1.$$
At this point the parity of $r$ comes into play.
If $r$ is even, then
$\delta(v_3^3,v_2^1)\le r$ and $\delta(v_3^3,v_2^2)\le r$ by
Lemma \ref{base_lemma},
thus for each $\ell\ne 3$ either
${\mathbf v}_3^3(\ell)={\mathbf v}_2^1(\ell)$
or ${\mathbf v}_3^3(\ell)={\mathbf v}_2^2(\ell)$ and
$\delta(v_3^3,v_2^1)=\delta(v_3^3,v_2^2)=r$. Accordingly
${\mathbf v}_3^3(1)=2,\; {\mathbf v}_3^3(2)=1$ and in the positions $\ell>3$
${\mathbf v}_3^3$ has $r\!-\!2$ ones but
${\mathbf v}_2^1$ has $r\!-\!3$ ones and
${\mathbf v}_2^2$ has $r\!-\!1$ ones,
leading to a
contradiction by Statement \ref{parity}.
If $r$ is odd, then we consider $v_3^{2r-1}$ instead of $v_3^{3}$,
here for each $\ell\ne 2r\!-\!1$ either
${\mathbf v}_3^3(\ell)={\mathbf v}_2^1(\ell)$
or ${\mathbf v}_3^3(\ell)={\mathbf v}_2^2(\ell)$, so
${\mathbf v}_3^3(1)=2,\; {\mathbf v}_3^3(2)=1$.
Now we focus to positions $\ell=3,\ldots,2r-2$ where
${\mathbf v}_2^{1},{\mathbf v}_2^{2},{\mathbf v}_3^{2r-1}$ have
$r\!-\!2$ ones.
By Statement \ref{parity} $\delta(v_3^{2r-1},v_2^{1})$
and $\delta(v_3^{2r-1},v_2^{2})$ is even, however, they should be exactly $r$
which is odd.
\ssqed

\begin{stm}\label{exact_dist}
Suppose $2\le i<i'\le r$.
Then
$$
\delta(v_i^j,v_{i'}^{j})=\left\{ \begin{array}{ll}
   r &\mbox{if $r$ is even}\\
   r-1 &\mbox{if $r$ is odd}\end{array}\right.
 $$
Also,  $\delta(v_i^j,v_{i'}^{j'})=r$ for any
$j\ne j'$. Moreover
$$
{\mathbf v}_i^{j}(j')=\left\{ \begin{array}{ll}
    {\mathbf v}_{i'}^{j'}(j)&\mbox{if $r$ is even}\\
    3-{\mathbf v}_{i'}^{j'}(j)&\mbox{if $r$ is odd}\end{array}\right.
$$
\end{stm}

\prf The first part is a consequence of Claim \ref{dist23} and
Statement \ref{parity}. If $\delta(v_i^j,v_{i'}^{j'})=r-1$ then by
Claim \ref{dist2comp} we get a vertex $a$ \wl in $V_{i}$ such that
$\delta(a,v_{i}^{j})=2r-1$, moreover ${\mathbf a}$ agrees with
${\mathbf  v}_{i'}^{j'}$ in all positions where
${\mathbf  v}_{i'}^{j'}$ and ${\mathbf  v}_{i}^{j}$ differ. Thus
${\mathbf a}(j')=3$, so we may call it $v_i^{j'}$ for getting a contradiction
by Statement \ref{notcompl_in_same_matrix}.
Having this, the last part is a consequence of Statement \ref{parity}.
\ssqed

We associate matrices to the selected vertices as follows.
For $i=2\ldots r$ let $A_i(j,j')=0$ if $j=j'$,
$A_i(j,j')=1$ if ${\mathbf  v}_{i}^{j}(j')=1$, and
$A_i(j,j')=-1$ if ${\mathbf  v}_{i}^{j}(j')=2$. These are
$(2r\!-\!1)\times(2r\!-\!1)$
matrices, we introduce the following operation for them.
$$A_i^*:=(-1)^r\cdot A_i^{\mathrm T}.$$
Now $A_{i'}=A_i^*$ for any $i\ne i'$ by the last part of Statement
\ref{exact_dist}. For the case $r\ge 4$ it is easy to complete the proof of
Claim \ref{nocomp}. Indeed, as $A_3=A_4=A_2^*$, we have e.g.,
${\mathbf  v}_{3}^{1}={\mathbf  v}_{4}^{1}$ contradicting to
Statement \ref{exact_dist}.
For the case $r=3$ we need some extra work. As $2r-1=5$ now, we have
$5\times 5$ matrices, and every row contains two $1$ and two $-1$ entries.
We define an auxiliary graph $G$ on
vertex set $[5]$. Let $ij$ is an edge iff $A_2(i,j)=A_2(j,i)=1$. We claim that
in this graph all five vertices have degree 1, leading to a contradiction.
If $d_G(i)=0$, then the $i$th row of $A_2$ is the negative of the $i$th
column. As $A_3=A_2^*$ we have
${\mathbf  v}_{2}^{i}={\mathbf  v}_{3}^{i}$ contradicting to
Statement \ref{exact_dist}.
If $d_G(i)\ge 2$, then the $i$th row of $A_2$ equals to the $i$th
column. As $A_3=A_2^*$ we have
$\delta({\mathbf  v}_{2}^{i},{\mathbf  v}_{3}^{i})=2r-2=4$, contradicting again to
Statement \ref{exact_dist}. Thus Claim \ref{nocomp} is proved.
\sqed

\bigskip

Combining the claims we conclude with the following corollary.

\begin{coroll}\label{cor} For any two vertices $v,w$ from different
classes, $\delta(v,w)=r$, and for any two vertices $u,v$
from the same class, $\delta(u,v)\le 2r-2$.
\end{coroll}

Now we are ready for finishing the proof of Theorem \ref{3plus2}. Select two
vertices $v\in V_1, \;w\in V_2$, \wl
${\mathbf v}=11\ldots 11, \; {\mathbf w}=22\ldots 2211\ldots1$ where
${\mathbf w}$ starts with exactly $r$ twos.
Accordingly, for a vector of length $2r\!-\!1$ we call its first part the
first $r$ coordinates, and its last part the last $r\!-\!1$ coordinates.
There exists a $y\in V_3$ with ${\mathbf y}(1)=3$, and let
$I=\{i \;|\; {\mathbf y}(i)=1\}$ and $\alpha=|I\cap [r]|$ (i.e., the number of
ones in its first part). Since $\delta(y,v)=r$, we have $|I|=r-1$, and since
$\delta(y,v)=\delta(y,w)$, ${\mathbf y}$ has exactly $\alpha$ twos in its
first part, so $\alpha\le \frac{r-1}{2}$.
There exists a
$J\subset [2r\!-\!1], \;J\supset I, \;|J|=r, \;|J\cap [r]| \le \frac{r-1}{2}$,
and by
the assumption ($r$ monochromatic components do not cover $V$)
there exist a vertex $v_J$ with the property
${\mathbf v_J}(j)\ne 1$ for all $j\in J$.

If $v_J\notin V_1\cup V_2$ then each of it coordinates
is 1 outside $J$, as $\delta(v_J,v)=r$. By the definition of $J$, it means that
${\mathbf v}_J$ has $\beta\ge \frac{r+1}{2}$ ones in the first part and
$\beta$ non-ones in the second part, so
$\delta(v_J,w)\ge \beta+\beta\ge r+1$, a contradiction.

Therefore $v_J\in V_1\cup V_2$, then $\delta(y,v_J)=r$
implies that ${\mathbf y}$ and
${\mathbf v_J}$ are equal outside $I$, with possibly one exception.
However, ${\mathbf y}(j)\ne 1$ for any $j\in [2r-1]-J$, consequently
${\mathbf  v_J}$ can have at most one coordinate that is 1. Thus
$v_J\in V_1$, and
$\delta(v_J,v)\ge 2r-2$, consequently by Corollary \ref{cor} it equals to
$2r-2$ and ${\mathbf v_J}$ has exactly one coordinate that is 1.

Let $I'=\{i\;|\; {\mathbf v_J}(i)={\mathbf w}(i)={\mathbf y}(i)\}$,
by Lemma \ref{gen_base_lemma} we have $|I'|\ge r-2$. However,
$I'\subset \{j\le r\;|\;   {\mathbf y}(j)=2\}$, and this latter set has
cardinality $\alpha\le \frac{r-1}{2}$, so $r-2\le \frac{r-1}{2}$,
i.e., $r\le 3$,
which leads to a contradiction, except for the case $r=3,\; \alpha=1$.

Now ${\mathbf v_J}(i)=1$ for an $i\le 3$, let $J'=\{i,4,5\}$ and define
$v_{J'}$ as ${\mathbf v_{J'}}(j)\ne 1$ if $j\in J'$. Now
$v_{J'}\in V_1\cup V_2$ because otherwise
$\delta(v_{J'},v)=\delta(v_{J'},w)=3$ would lead to a contradiction.
If $v_{J'}\in V_2$, then  ${\mathbf v_{J'}}(j)=1$ for $j\not\in J'$, now we
choose $z\in V_3$ with ${\mathbf z}(i)=1$. As $\delta(z,v)=3$, $\mathbf z$ has
exactly one other 1 but if its position is in $\{1,2,3\}$, then
$\delta(z,w)\ge 4$, and if in  $\{4,5\}$, then  $\delta(z,v_{J'})\ge 4$.

So $v_{J'}\in V_1$, consequently, by $\delta(v_J,w)=\delta(v_{J'},w)=3$, both
$v_J$ and $v_{J'}$ have 2 twos in the first part, let $\ell\in[3]$
the position of a common 2 and we now choose $t\in V_3$ with
${\mathbf t}(\ell)=3$.  Since $\delta(t,v)=\delta(t,w)=3$, $t$ has one 1 and one
2 in the first part, and one 2 in the second part, contradicting to
$\delta(v_J,t)=\delta(v_{J'},t)=3$.
\qed

\begin{coroll}\label{coroll}
If $r\ge 3$, then $\cov(r,k)=1$ for every $1\le k\le r$,\\
  $\cov(r,k)=k-r+1$ for every  $r\le k\le 2r-1$,\\
  and for any $k\ge 2r$ we have $k-r+1\le \cov(r,k)\le k-r+2$.
\end{coroll}

\section{Generalized complete uniform hypergraphs}\label{trans}

\begin{defi}
A  hypergraph is called $(r,\ell)$-{\bf partite} if the
ground set $V$ is partitioned into nonempty classes
$V_1\cup\ldots\cup V_r$, and
no edge intersects any $V_i$ in more than $\ell$ vertices.
A hypergraph is complete $r$-uniform  $(r,\ell)$-partite if its edge set
consists of all $r$-tuples intersecting each class in at most $\ell$ vertices.
An edge of an $(r,\ell)$-partite  hypergraph is called \emph{friendly} if it
intersects  at most one class in exactly $\ell$ vertices; otherwise we call it
\emph{unfriendly}.
An $r$-uniform $(r,\ell)$-partite  hypergraph is called
{\bf semicomplete} if its edge set consists of all $r$-tuples
intersecting at most one class in exactly $\ell$ vertices (that is, it consists
of the friendly edges of the complete $r$-uniform
$(r,\ell)$-partite hypergraph).
An $r$-uniform $(r,\ell)$-partite  hypergraph is called
{\bf rich} if it contains all edges of the semicomplete hypergraph.
\end{defi}

Among  $r$-uniform hypergraphs
the complete $(r,1)$-partite hypergraphs are the complete
$r$-partite ones and complete $(r,r)$-partite hypergraphs are the
complete ones. The complete $(r,r\!-\!1)$-partite hypergraphs are
also interesting, containing all $r$-tuples of $V$ except those that are
contained in some $V_i$.  The purpose of this section is to build a bridge
between the two known extreme cases ($\ell=r$ was solved in \cite{KI},
$\ell=1$ was handled in the previous section).

For $1\le \ell \le r$, let $\cov(r,\ell,k)$ denote the minimum number of
monochromatic components
needed to cover the vertex set of any complete $r$-uniform $(r,\ell)$-partite
hypergraph in any spanning $k$-coloring.

\begin{conj}\label{conj3}
  $$\cov(r,\ell,k)= 1+\Bigl\lfloor \frac{k-r+\ell-1}{\ell}\Bigr\rfloor$$
  for every $r\ge 3,\; k\ge 1+r-\ell, \; 1\le \ell\le r$.
\end{conj}

We start with giving the lower bound.

\begin{constr}\label{rl_constr}
  This construction is a straightforward generalization of Construction
  \ref{basic_constr}. We have $r,k,\ell$ fixed with $k\ge r+1\ge 4$ and
  $1\le \ell \le r$, and let $q= \lfloor\frac{k-r+\ell-1}{\ell} \rfloor$ and
  $k'=q\cdot \ell +r-\ell+1\le k$.
  First we fix the sizes and labels of the classes.
  $|V_1|={k'\choose q}$ and elements $V_1$ are labeled with the
  $q$-element subsets of $[k']$.
  For $2\le j\le r$ set $V_j$ is a disjoint union of
  $V_j=A_j^1\cup \dots \cup A_j^{k'}$ where $|A_j^i|={k'-1\choose q-1}$, all
  elements of $A_j^i$ are labeled with set $\{i\}$ and have an arbitrary fixed
  linear order. Now take an arbitrary rich $r$-uniform
  $(r,\ell)$-partite  hypergraph
  $\KR$ on $V=V_1\cup\ldots\cup V_r$, we are going to define
  a spanning $k'$-coloring of its edges.

First we define special edges of color $i$ for any $i\in [k']$. Consider the
set $W_i$ of ${k'-1\choose  q-1}$
vertices of $V_1$ associated to $q$-sets of $[k']$ containing $i$.

\emph{Special edges of color $i$} are the
${k'-1\choose  q-1}$ edges whose vertex from $W_i$ is the $\ell$-th in
lexicographic order, and for all
$2\le j\le r$ whose vertex from $V_j$ is the $\ell$-th
in the fixed linear order of $A_j^i$ for $\ell=1,\dots,{k'-1\choose q-1}$.
Thus special edges of color $i$ form a matching for all $i$.

\emph{Non-special edges} with vertices $v_1,\dots,v_r$ get
  their color as
  the smallest $c\in [k']$ such that $c$ is not in the
  union of sets associated to $v_1,\dots,v_r$.

  Note that every non-special
$r$-tuple $v_1,\dots,v_r$ gets a color because the conditions forbid at most
  $\ell\cdot q+(r-\ell)<k'$ colors.
  Observe also that a special edge of color $i$ is always
disjoint from any other edge of color $i$. Consequently a special edge of
color $i$ forms a monochromatic component of color $i$ having $r$ vertices,
we call them small monochromatic components.

We claim that the coloring is spanning.
Suppose first that $v\in V_1$ representing the set
$Q_1\subset [k']$.  For any $i\in Q_1$,  $v$
is in a special edge of color $i$. On the other hand, for any $i\not\in Q_1$ we
can select vertices $v_2,\ldots,v_\ell\in V_1$ with associated $q$-sets
$Q_2,\ldots,Q_\ell\subseteq [k']-\{i\}$,
such that for every $j\ne j'$ sets $Q_j$ and $Q_{j'}$ are
disjoint. Then we may select $v_{\ell+1},\ldots,v_r$ from
$V_2,\ldots,V_{r-\ell+1}$, such that the associated one-element subsets are
distinct, and are subsets of $[k']-\{i\}-\cup Q_j$.
Now the union of the associated
sets of our selected $r$-tuple is $[k']-\{i\}$, thus it was colored by $i$.

On the other hand, let $v_r\in A_j^i$
for some $2\le j\le r,\; 1\le i \le k'$. Clearly
$v_r$ is in a special edge of color $i$.
For any  $1\le c \le k'$ if $c\ne i$, then we can take
 vertices $v_1,\ldots,v_\ell\in V_1$ with associated $q$-sets
 $Q_1,\ldots,Q_\ell\subseteq [k']-\{i\}-\{c\}$,
 such that for every $j\ne j'$ sets $Q_j$ and $Q_{j'}$ are
disjoint. Then we may select $v_{\ell+1},\ldots,v_{r-1}$ from
$V_2\cup\ldots\cup V_r-V_j$, such that the associated one-element subsets are
distinct, and are subsets of $[k']-\{i\}-\{c\}-\cup Q_j$.
Now the union of the associated
sets of our selected $r$-tuple is $[k']-\{c\}$, thus it was colored by $c$.
\end{constr}

\begin{thm}\label{r_l_l}
  $\cov(r,\ell,k)\ge 1+\lfloor \frac{k-r+\ell-1}{\ell}\rfloor$
  for every $r\ge 3,\; k\ge 1+r-\ell, \; 1\le \ell\le r$.
\end{thm}

\prf The statement is obvious if $k\le r$.
Consider Construction \ref{rl_constr}.
Note that the union of at most $q=\lfloor\frac{k-r+\ell-1}{\ell} \rfloor$
large monochromatic components
do not cover $V_1$. Let their colors are $c_1,\ldots,c_s$ with $s\le q$, and
take any $q$-set that contains $\{c_1,\ldots,c_s\}$; the vertex in $V_1$
associated to this set is not covered.

The uncovered vertices of $V_1$ must be covered by small monochromatic
components, and every such component can contain just one vertex of $V_1$.
Therefore we need ${k'-s\choose q-s}>q-s$ small monochromatic components to
cover them, thus altogether we need more than $s+(q-s)=q$ monochromatic
components to cover all vertices.
\qed

\begin{remark}
The basic idea of the above construction is from \cite{KI} where the
constructed coloring for complete $r$-uniform hypergraphs is not spanning
(this was not an issue of that paper). Here, when $\ell=r$, we gave another
construction for complete $r$-uniform hypergraphs where the coloring is spanning.
\end{remark}

\begin{thm}\label{r_l_u}
  $\cov(r,\ell,k)\le 1+\lfloor\frac{k-r+\ell-1}{\ell}\rfloor$
  for every $r\ge 3,\; k\ge 1+r-\ell, \; 2\le \ell\le r$.
\end{thm}

\prf
The proof goes similarly as in the proof of Theorem \ref{lemma_essential}.
Fix the nonempty classes $V_1,\ldots,V_r$ and take any rich
$r$-uniform $(r,\ell)$-partite
hypergraph $\KR$ with a spanning $k$-coloring of its edges.
We are going to show by induction on $k$ that
$c(\KR)\le 1+\lfloor\frac{k-r+\ell-1}{\ell}\rfloor$.
The cases $k\le r$ are obvious.

Let $e=\{u_1,\ldots,u_r\}$ be an essential edge of $\KR$
colored by 1, if no such edge exists, then recolor edges having
color 1 and use induction. Until there exists an essential friendly edge
colored by 1, we choose that edge for $e$. If all essential edges
colored by 1 are  unfriendly, then simply delete them  from
$\KR$ getting a $(k\!-\!1)$-colored rich hypergraph, where the coloring is
still spanning, so we are done by induction.

So $e$ is a friendly essential edge, \wl $\ell\ge|e\cap V_1|\ge |e\cap V_j|$
for all $j$. As $e$ is friendly, we also have $|e\cap V_j|< \ell$
for $j>1$.
Take $R_{u_1},\ldots,R_{u_r}$, where $R_{u_j}=e-\{u_j\}$,
for any $i\ne j$ we have $\Col(R_{u_i})\cap\Col(R_{u_j})=\{1\}$,
so there is a $j$ with
$|\Col(R_{u_j})|\le 1+\lfloor \frac{k-1}{r}\rfloor$.

First consider the case $|R_{u_j}\cap V_1|<\ell$
(note that this is always true for $\ell=r$).
We also emphasize here that for this case we do not need
the coloring to be spanning.
For any vertex $v\in V$ the set $R_{u_j}\cup\{v\}$
is a friendly edge of $\KR$, consequently
the monochromatic components of colors in $\Col(R_{u_j})$ containing
$R_{u_j}$ cover the whole $V$. We need to prove
$\lfloor \frac{k-1}{r}\rfloor\le\lfloor\frac{k-r+\ell-1}{\ell}\rfloor$.
For $k-1<r$ both are zero, otherwise $(r-\ell)(k-1)\ge (r-\ell)r$, so
$\frac{k-1}{r}\le \frac{k-r+\ell-1}{\ell}$.

So we are left with the case $|R_{u_j}\cap V_1|= \ell$.
There are two possibilities.
Either one of $\Col(R_{u_i})=\{1\}$ for an $i>1$,
in this case the monochromatic component
containing $u_i$ and colored by 1 covers $V$ because it covers $V-V_1$,
as for all $v\in V-V_1$ the set $e-\{u_i\}\cup\{v\}$ is an edge of $\KR$,
and (using that the coloring is spanning),
every $w\in V_1$ is incident to an edge colored by 1 and this edge meets
$V-V_1$.

Otherwise $|\Col(R_{u_i})|\ge 2$ for all $i>1$, so by the pigeonhole principle
there is a $2\le i\le\ell$ with
$|\Col(R_{u_i})|\le 1+\lfloor \frac{k-1-(r-\ell)}{\ell}\rfloor$, and
the monochromatic components of colors in $\Col(R_{u_i})$ containing
$e-\{u_i\}$ cover the whole $V$ because $e-\{u_i\}\cup\{v\}$ is an edge of
$\KR$ for every $v\in V-e$.
\qed

Summarizing the results of this section and Corollary
\ref{coroll}, we proved Conjecture \ref{conj3} for almost all cases.
We also proved that Conjecture \ref{conj3} is equivalent to Conjecture
\ref{conj}.

\begin{thm}[Main theorem]\label{r_l_e}
  $$\cov(r,\ell,k)= 1+\Bigl\lfloor \frac{k-r+\ell-1}{\ell}\Bigr\rfloor$$
  for every $r\ge 3,\; k\ge 1+r-\ell, \; 1\le \ell\le r$, \emph{except} when
  $\ell=1$ and $k\ge 2r$, where only $1+\lfloor \frac{k-r+\ell-1}{\ell}\rfloor$
  $\le \cov(r,\ell,k)\le 2+\lfloor \frac{k-r+\ell-1}{\ell}\rfloor$ was proved.
\end{thm}

\section{Open problems}

Besides the missing case $(k\ge 2r)$ of Conjecture \ref{conj} and the above
mentioned conjecture of Gy\'arf\'as and Lehel (stating that $\COV(2,k)=2k-2$),
we list some more open problems.

In \cite{CFGYLT} it is shown that $2k-2\le \COV(2,k)\le 2k-1$. Much less is
known about $\cov(2,k)$. The best known upper bound is still $2k-1$ but no
reasonable lower bound is known. The second author conjectures that
$2k-4\sqrt{k}\le \cov(2,k)$.

For $r\ge 3$
we did not study $\COV(r,\ell,k)$ (that is similar to $\cov(r,\ell,k)$ but the
coloring need not to be a spanning one), it was only determined for $\ell=1$
(see Theorem \ref{nonspanning}) and for $\ell=r$
(either in \cite{KI} or in the proof of Theorem \ref{r_l_u}).

\bigskip

We can naturally generalize further.

\begin{defi}
  For $1\le \ell \le r\le R\ell$,
  let $\cov(r,R,\ell,k)$ denote the minimum number of
monochromatic components
needed to cover the vertex set of any complete $r$-uniform $(R,\ell)$-partite
hypergraph in any spanning $k$-coloring.
\end{defi}

Determining $\cov(r,R,\ell,k)$ for all possible ranges seems to be very
challenging. At the moment we do not have a conjecture about the value of
$\cov(2,3,1,k)$.


\begin{thebibliography}{99}

\bibitem{GY1} A. Gy\'arf\'as, Partition covers and blocking sets in
  hypergraphs, {\em MTA SZTAKI tanulm\'anyok} {\bf 71} (1977) (in Hungarian)
\bibitem{CFGYLT} G. Chen, S. Fujita, A. Gy\'arf\'as, J. Lehel, \'A. T\'oth,
  Around a biclique cover conjecture, {\em arxiv:1212.6861}
\bibitem{EGYP} P. Erd\H os, A. Gy\'arf\'as, L. Pyber, Vertex coverings by
  monochromatic cycles and trees, {\em Journal of Combinatorial Theory B}
 {\bf 51.} (1991) 90-95.
\bibitem{FFGYT1} S. Fujita, M. Furuya, A. Gy\'arf\'as, \'A. T\'oth, Partition
  of graphs and hypergraphs into monochromatic connected parts,
{\em Electronic Journal of Combinatorics} {\bf 19} P27.

\bibitem{FFGYT2} S. Fujita, M. Furuya, A. Gy\'arf\'as, \'A. T\'oth, A note on
  covering edge colored hypergraphs by monochromatic components,
{\em Electronic Journal of Combinatorics} (2014) {\bf 21} P33.

\bibitem{GY2} A. Gy\'arf\'as, Vertex covers by monochromatic pieces - A survey
  of results and problems, {\em Discrete Mathematics}, to appear (2016)

\bibitem{HE} J. R. Henderson, Permutation Decomposition of (0-1)-Matrices and
Decomposition Transversals, Ph.D. thesis, Caltech, 1971.


\bibitem{KI} Z. Kir\'aly, Monochromatic components in edge-colored complete
  hypergraphs, {\em European Journal of Combinatorics} {\bf 35} (2013)
  374-376.

\end{thebibliography}
\end{document}